# Nonexistence of smooth Levi-flat hypersurfaces in complex projective spaces of dimension $\geq 3$

By Yum-Tong Siu*

In this paper we prove the following theorem.

MAIN THEOREM. *Let $n \geq 3$ and $m \geq \frac{3n}{2} + 7$. Then there exists no $C^m$ Levi-flat real hypersurface $M$ in $\mathbf{P}_n$.*

The condition that $M$ is Levi-flat means that when $M$ is locally defined by the vanishing of a $C^m$ real-valued function $f$, at every point of $M$ the restriction of $\partial\bar\partial f$ to the complex tangent space of $M$ is identically zero.

The case of the nonexistence of $C^\infty$ Levi-flat real hypersurface in $\mathbf{P}_2$ is motivated by problems in dynamical systems in $\mathbf{P}_2$ (see [LN]). For some technical reason the method of this paper at this point can only yield the result for the case of $\mathbf{P}_n$ with $n \geq 3$. That technical reason will be explained later in this introduction.

By slicing we can reduce the case of a general $n$ to the case of $n = 3$ with a weaker assumption of the order of differentiability of the Levi-flat, real hypersurface. The same proof works when $\mathbf{P}_n$ is replaced by an irreducible compact Hermitian symmetric manifold $X$ of complex dimension $n$ whose bisectional curvature is $(n-2)$-nondegenerate. (The definition for the (strong) nondegeneracy of the bisectional curvature is given in Definition 2.3.)

THEOREM 1. *Let $m \geq \frac{3n}{2} + 7$. Then there exists no $C^m$ Levi-flat, real hypersurface $M$ in an irreducible compact Hermitian symmetric manifold $X$ of complex dimension $n$ whose bisectional curvature is $(n-2)$-nondegenerate.*

For the special case when $M$ is assumed to be real-analytic, the Main Theorem was proved by Lins Neto [LN]. In that proof the real-analyticity is required to conclude the extension in an appropriate way of the structure of $M$ to a neighborhood. The proof of the differentiable case here requires a completely different approach and is reduced to the regularity problem in solving the $\bar\partial_b$ equation for a $(0,1)$-form on the Levi-flat hypersurface $M$. Such a regularity statement is not expected to hold for a general Levi-flat manifold.

*Partially supported by a grant from the National Science Foundation.

The main idea of the proof of Theorem 1 is to get a contradiction from the complex normal bundle $T_X^{1,0}/T_M^{1,0}$ of the $C^m$ Levi-flat real hypersurface $M$ in $X$. The contradiction comes from the following two facts. The first fact is that, when restricted to any local holomorphic leaf of the foliation of $M$, the complex normal bundle $T_X^{1,0}/T_M^{1,0}$ as a quotient of the tangent bundle $T_X^{1,0}$ carries some positivity. The second fact is that the complex normal bundle $T_X^{1,0}/T_M^{1,0}$ of the real hypersurface $M$ is topologically trivial and admits a $C^2$ Hermitian metric along its fibers which has zero curvature when restricted to any local holomorphic leaf of the foliation of $M$. The two facts together give a Hermitian metric of the trivial complex line bundle over $M$ whose restriction to any local holomorphic leaf of the holomorphic foliation of $M$ carries some positivity. A contradiction occurs when one considers the point where $-\log$ of that Hermitian metric of the trivial line bundle achieves its maximum.

The main difficulty of the proof is the second fact, which is the analog, for the Levi-flat hypersurface $M$, of Kodaira's result that a $d$-exact $(1,1)$-form on a compact Kähler manifold is the $\partial\bar\partial$ of a function. The proof of such an analog hinges on the regularity problem in solving the $\bar\partial_b$ equation for a $C^{m-1}$ $(0,1)$-form $\alpha$ on $M$. The method used in our solution of the regularity problem in our case requires the dimension $n$ of $\mathbf{P}_n$ to be at least 3, though the second fact is expected to hold even for the case of $n = 2$.

Our method to prove the second fact is as follows. First $\alpha$ is extended to a $C^{m-1}$ $(0,1)$-form $\tilde\alpha$ on $X$ whose $\bar\partial$ vanishes to order $m-2$ at $M$. Then the $L^2$ estimates of $\bar\partial$ on $X - M$ give a $C^{m-1}$ $(0,1)$-form $\beta$ on $X - M$ so that $\bar\partial\beta = \bar\partial\tilde\alpha$ on $X$ in the sense of currents and the quotient of $\beta$ by the $(m-2)$-th power of the distance function to $M$ is $L^2$ in a neighborhood of $M$. This step involves use of the curvature property of $X$ to give $X - M$ a suitable complete Kähler metric constructed from the distance function to $M$ so that $-\log$ of the distance function to $M$ satisfies the condition that the sum of at least $n-2$ of the eigenvalues of its complex Hessian with respect to the Kähler metric is bounded from below by a positive constant near $M$. It is in this step the condition that the $n$ in $\mathbf{P}_n$ is at least 3 is needed. This step corresponds to the $L^2$ analog of the vanishing of the second cohomology with compact support for Stein manifolds of complex dimension at least 3. After this step one solves for $\gamma$ in $\bar\partial\gamma = \tilde\alpha - \beta$ on $X$ and shows that the restriction of $\gamma$ to $M$ is $C^p$ for $p \leq n - \frac{3n}{2} - 5$ and $\bar\partial_b\gamma = \alpha$. The conclusion that the restriction of $\gamma$ to $M$ is $C^p$ is derived by use of the pole order of the explicit kernel for the local solution of the $\bar\partial$ equation and the fact that the restriction to a real hypersurface of the solution of a $\bar\partial$ is $L^1$ when the right-hand side is $L^2$.

The following is a more intuitive (but technically imprecise) way to explain why our proof of the Main Theorem requires $n \geq 3$. If the $\bar\partial_b$-equation can be solved with regularity on $M$, then type considerations yield readily the second fact above, which says that any $d$-exact 2-form on $M$ of type $(1,1)$ is $\partial_b \bar{\partial}_b$-exact with regularity. The solvability of the $\bar{\partial}_b$-equation with regularity on $M$ is handled by the $L^2$ analog of the following exact sequence:

$$H^1(X, \mathcal{O}_X) \to H^1(M, \mathcal{O}_X|M) \to H^2_{\text{compact}}(X - M, \mathcal{O}_X),$$

where $H^2_{\text{compact}}(X - M, \mathcal{O}_X)$ is the cohomology group with compact support. The solvability of the $\bar{\partial}_b$-equation with regularity on $M$ is analogous to the vanishing of $H^1(M, \mathcal{O}_X|M)$. The cohomology group $H^1(X, \mathcal{O}_X)$ is zero due, for example, to the simple connectedness of $X$. The vanishing of $H^2_{\text{compact}}(X - M, \mathcal{O}_X)$ and its $L^2$ analog require $n \geq 3$ when $X = \mathbf{P}_n$.

In the general case of the nonvanishing of the $(0,1)$-cohomology group, the proof of Kodaira's lemma on $\partial\bar{\partial}$-exactness uses the complex conjugation relation between the $(1,0)$ and $(0,1)$ cohomology groups in Hodge theory to overcome the difficulty. The difficulty of our case of $\mathbf{P}_2$ with the nonvanishing of $H^2(\mathbf{P}_2 - M, \mathcal{O}_{\mathbf{P}_2})$ (with $L^2$ bound) should probably be overcome with a technique corresponding to the complex conjugation property in Hodge theory used in the proof of Kodaira's lemma.

Theorem 1 should be generalizable to the nonexistence of any smooth, Levi-flat, real submanifold of real codimension $q$ and constant complex dimension $n-q$ in an irreducible compact Hermitian symmetric manifold of complex dimension $n$ whose bisectional curvature is strongly $(n-q-1)$-nondegenerate (possibly with the additional assumption that the determinant line bundle of the complex normal bundle of the submanifold is topologically trivial). A good modification of the last step in the argument to solve with regularity the $\bar{\partial}_b$ equation on the Levi-flat real submanifold is needed, because the restriction to a real submanifold of higher real codimension of the solution of the $\bar{\partial}$-equation is no longer $L^1$ when the right-hand side is $L^2$. A couple of other easier modifications are also needed for such a generalization and they will be discussed in Section 8 below.

For the conjectured case of the Main Theorem for $\mathbf{P}_2$, Ohsawa and Sibony [O-S] introduced a method to get a contradiction by constructing an infinite number of linearly independent holomorphic sections of a suitable line bundle $L$ over $\mathbf{P}_2$ by using a suitable complex line $H$ in $\mathbf{P}_2$ and extending smooth sections of $L$ over $H \cap M$ first to $M$ and then to all of $\mathbf{P}_2$. Their method reduces the problem to the regularity problem of solving the $\bar{\partial}_b$-equation for an $L$-valued $(0,1)$-form on $M$, which still remains open. Their construction of an infinite number of linearly independent holomorphic sections of a suitable line bundle $L$ over $\mathbf{P}_2$ depends heavily on the fact that $H \cap M$ is a real curve and hence requires the dimension $n$ of $\mathbf{P}_n$ to be 2.

The technique given in this paper introduces a new way of obtaining regularity in solving the $\bar{\partial}_b$ equation on a Levi-flat hypersurface. Unfortunately the conditions required for its use are very restrictive. As a matter of fact, in the setup of this paper there cannot be any Levi-flat hypersurface which sat-isfies such very restrictive conditions. It is worth exploring to see what cases of weakly pseudoconvex hypersurfaces this new way of solving the $\bar\partial_b$ equation with regularity is applicable to.

Note that if an invariant transversal measure is assumed for the foliation of $M$, the usual techniques of Hodge theory apply to $M$ to yield readily the $\partial_b\bar\partial_b$-exactness of a $d$-exact $(1,1)$-form on $M$ with regularity.

The rest of the paper is devoted to the proof of Theorem 1.

## 1. The $L^2$ vanishing theorem for negative line bundles on complete Kähler manifolds

Both the statement and the proof of the $L^2$ vanishing theorem in this section are from known standard techniques and they are presented here just to get the result precisely in the form needed here.

THEOREM 1.1 (The $L^2$ vanishing theorem for negative line bundles on complete Kähler manifolds). *Let $q$ be a positive integer. Let $X$ be a complete Kähler manifold of complex dimension $n > q$ and $L$ be a holomorphic line bundle over $X$ with a $C^m$ Hermitian metric along its fibers for some $m \geq 2$. Assume that $X$ is Stein or satisfies the weaker condition that every compact subset $K$ of $X$ admits an open neighborhood $U$ such that the global holomorphic vector fields on $U$ generate the tangent space of $X$ at every point of $U$. Let $\kappa_0$ be a positive number. Further assume that at every point of $X$ the sum of any $n-q$ eigenvalues of the curvature form of $L$ with respect to the Kähler metric of $X$ is bounded from above by $-\kappa_0$. Then for any measurable $\bar\partial$-closed $L$-valued $(0,q)$-form $\omega$ on $X$ with $\int_X |\omega|^2$ finite, there exists a measurable $L$-valued $(0, q-1)$-form $u$ on $X$ such that $\bar\partial u = \omega$ in the sense of currents and $\int_X |u|^2 \leq \frac{1}{\kappa_0} \int_X |\omega|^2$. Moreover, the solution $u$ is $C^{k+1}$ on any open subset of $X$ where $\omega$ is $C^k$ if $k+1 \leq m$.*

*Proof.* Let $g_{\alpha\bar\beta}$ be the complete Kähler metric of $X$ and let $h$ be the Hermitian metric of $L$ and $\Theta_{\alpha\bar\beta} = -\partial\bar\partial \log h$ be the curvature form of $h$. Here and in the rest of this paper we use the summation of convention of summing over a symbol appearing at the same time as a superscript and a subscript, without writing down the summation sign. Let $\operatorname{Tr}\Theta = \Theta_{\alpha\bar\beta} g^{\alpha\bar\beta}$. For a $C^\infty$ $L$-valued $(0,q)$-form $\varphi$ on $X$ we define the operator $\varphi \mapsto \tilde\Theta\varphi$ by

$$\left(\tilde\Theta\varphi\right)_{\bar j_1\cdots\bar j_q} = \sum_{\nu=1}^{q} \Theta^{\bar\ell}{}_{\bar j_\nu} \varphi_{\bar j_1\cdots(\bar\ell)_\nu\cdots\bar j_q} - (\operatorname{Tr}\Theta)\varphi_{\bar j_1\cdots\bar j_q},$$

where the subscript $(\bar\ell)_\nu$ in the first term on the right-hand side means that the index in the $\nu$-th position is replaced by the index $\bar\ell$. At an arbitrarily prescribed point $P$ in $X$ we can choose local coordinates at $P$ so that $g_{\alpha\bar\beta}$ is equal to the Kronecker delta $\delta_{\alpha\beta}$ and the $(1,1)$-form $\Theta_{\alpha\bar\beta}$ is in diagonal form with $\Theta_{\alpha\bar\beta} = \lambda_\alpha \delta_{\alpha\beta}$, where the $\lambda_\alpha$ are the eigenvalues of $\Theta_{\alpha\bar\beta}$ with respect to $g_{\alpha\bar\beta}$. Since $\sum_{k=1}^{n-q} \lambda_{\nu_k} \leq -\kappa_0$ for all $1 \leq \nu_1 < \cdots < \nu_{n-q} \leq n$, it follows that $\left\langle \tilde\Theta \varphi, \varphi \right\rangle \geq \kappa_0 |\varphi|^2$ at the point $P$, where $\langle \cdot, \cdot \rangle$ and $|\cdot|$ denote respectively the pointwise inner product and the pointwise norm. By direct computation,

$$(\Box \varphi)_{\bar j_1 \cdots \bar j_q} = -g^{i\bar j}\nabla_{\bar j}\nabla_i \varphi_{\bar j_1 \cdots \bar j_q} + \left(\tilde\Theta \varphi\right)_{\bar j_1 \cdots \bar j_q}$$

(see e.g., [Siu2, (1.3.4)] where the curvature is defined with a different sign convention). When $\varphi$ has compact support in $X$, integration by parts yields

$$\|\bar\partial \varphi\|_X^2 + \|\bar\partial^* \varphi\|_X^2 = \|\nabla \varphi\|_X^2 + \left(\tilde\Theta \varphi, \varphi\right)_X \geq \kappa_0 \|\varphi\|_X^2,$$

where $(\cdot,\cdot)_X$ and $\|\cdot\|_X$ denote respectively the global inner product and the global norm over $X$. Let $\psi$ be any $L$-valued $(0,q)$-form with compact support $K$ in $X$ which is $L^2$ and which belongs to the domain of $\bar\partial$ and the domain of $\bar\partial^*$. Since $K$ admits an open neighborhood $U$ so that global holomorphic vector fields of $U$ generate the tangent space at every point of $U$, we can use biholomorphisms of a neighborhood of $K$ defined by global holomorphic vector fields of $U$ to smooth out $\psi$. By the argument of Friedrichs's lemma [Frie] we conclude that there exists a sequence of $C^\infty$, $L$-valued $(0,q)$-forms $\varphi_\nu$ with common compact support in $X$ which approach $\psi$ in the graph norm in the sense that

$$\|\varphi_\nu - \psi\|_X^2 + \|\bar\partial(\varphi_\nu - \psi)\|_X^2 + \|\bar\partial^*(\varphi_\nu - \psi)\|_X^2$$

approaches $0$ as $\nu \to \infty$. Hence

$$\|\bar\partial \psi\|_X^2 + \|\bar\partial^* \psi\|_X^2 \geq \kappa_0 \|\psi\|_X^2.$$

We now remove the condition that $\psi$ has compact support, but continue to assume that $\psi$, $\bar\partial \psi$, and $\bar\partial^* \psi$ are all $L^2$ on $X$. Since the Kähler manifold $X$ is complete, given any compact subset $K$ in $X$ we can find a $C^\infty$ function $0 \leq \rho \leq 1$ with compact support in $X$ which is identically $1$ in a neighborhood of $K$ such that $|d\rho| \leq 1$ at every point of $X$. Then from

$$\bar\partial(\psi - \rho\psi) = (1-\rho)\bar\partial \psi - (\bar\partial \rho)\psi$$

it follows that

$$\|\bar\partial(\psi - \rho\psi)\|_X \leq \|\bar\partial \psi\|_{X-K} + \|\psi\|_{X-K}.$$

From

$$\bar\partial^*(\psi - \rho\psi) = (1-\rho)\bar\partial^* \psi - \bar\partial \rho \vdash \psi,$$

where $\vdash$ denotes the interior product, it follows that

$$\|\bar\partial^*(\psi - \rho\psi)\|_X \leq \|\bar\partial^* \psi\|_{X-K} + n\|\psi\|_{X-K}.$$

As $K$ goes through an increasing sequence of compact subsets which exhausts $X$, we conclude that the inequality

(1.1.1) $$\|\bar\partial\psi\|_X^2 + \|\bar\partial^*\psi\|_X^2 \geq \kappa_0 \|\psi\|_X^2$$

is still valid for all $L^2$ $\psi$ with both $\bar\partial\psi$ and $\bar\partial^*\psi$ also $L^2$. By writing $\psi = \psi_1 + \psi_2$ with $\psi_1 \in \operatorname{Ker}\bar\partial$ and $\psi_2 \in (\operatorname{Ker}\bar\partial)^\perp \subset \operatorname{Ker}\bar\partial^*$, we have

$$\begin{aligned} |(\psi,\omega)_X| &= |(\psi_1,\omega)_X| \leq \|\psi_1\|_X \|\omega\|_X \\ &\leq \frac{\|\omega\|_X}{\sqrt{\kappa_0}} \|\bar\partial^*\psi_1\|_X = \frac{\|\omega\|_X}{\sqrt{\kappa_0}} \|\bar\partial^*\psi\|_X \end{aligned}$$

for all $L^2$ $\psi$ with both $\bar\partial\psi$ and $\bar\partial^*\psi$ also $L^2$. By Riesz's representation theorem applied to the functional

$$\bar\partial^*\psi \mapsto (\psi,\omega),$$

we conclude that there exists uniquely an $L^2$, $L$-valued $(0,q-1)$-form $u$ on $X$ perpendicular to $\operatorname{Ker}\bar\partial$ such that $(\bar\partial^*\psi, u)_X = (\psi,\omega)_X$ for all $L^2$ $\psi$ with both $\bar\partial\psi$ and $\bar\partial^*\psi$ also $L^2$. Hence $\bar\partial u = \omega$. Moreover, $\|u\|_X \leq \frac{\|\omega\|_X}{\sqrt{\kappa_0}}$.

To get the final statement of $u$ being $C^{k+1}$ wherever $\omega$ is $C^k$, we need to use the ellipticity of $\bar\partial\bar\partial^* + \bar\partial^*\bar\partial$. From (1.1.1) we obtain

$$\|(\bar\partial\bar\partial^* + \bar\partial^*\bar\partial)^{-1}\varphi\|_X \leq \frac{1}{\sqrt{\kappa_0}}\|\varphi\|_X,$$

which implies that there exists $v$ such that $\bar\partial v \in \operatorname{Dom}\bar\partial^*$ and $\bar\partial^*v \in \operatorname{Dom}\bar\partial$ and $(\bar\partial\bar\partial^* + \bar\partial^*\bar\partial)v = \omega$ (cf., [F-K, (1.3.8)]). By applying $\bar\partial$ to the last equation and taking the inner product with $\bar\partial v$, we conclude that $\bar\partial^*\bar\partial v = 0$ and $\bar\partial\bar\partial^*v = \omega$. Since $\bar\partial^*v$ is orthogonal to $\operatorname{Ker}\bar\partial$, it follows that $u = \bar\partial^*v$. From the ellipticity of the operator $\bar\partial\bar\partial^* + \bar\partial^*\bar\partial$ it follows that $u$ is $C^{k+1}$ wherever $\omega$ is $C^k$ if $k+1 \leq m$. □

## 2. Lower bound of the suhbarmonicity of the distance function to a Levi-flat hypersurface

We start with the following simple standard lemma.

LEMMA 2.1 (The integral form of subharmonicity). *For a function $f$ which is $C^2$ in a neighborhood of $0$ in $\mathbf{C}$,*

$$(\Delta f)(0) = \lim_{r \to 0} \frac{4}{r^2}\left(\frac{1}{2\pi}\int_{\theta=0}^{2\pi} f(re^{\sqrt{-1}\theta})d\theta - f(0)\right).$$

*Proof.* Let $B(r)$ denote the open disk of radius $r$ in $\mathbf{C}$ centered at the origin and let $\frac{\partial}{\partial \vec{n}}$ denote the differentiation in the direction of the unit outward normal. When we apply $\int_{\rho=0}^{r} \frac{d\rho}{\rho}$ to both sides of

$$\int_{B(\rho)} \Delta f = \int_{\partial B(\rho)} \left(\frac{\partial}{\partial \vec{n}} f\right) \rho d\theta,$$

it follows that

$$\frac{1}{2\pi} \int_{\rho=0}^{r} \frac{d\rho}{\rho} \int_{B(\rho)} \Delta f = \frac{1}{2\pi} \int_{\theta=0}^{2\pi} f(re^{\sqrt{-1}\theta}) d\theta - f(0).$$

For any positive number $\varepsilon$

$$((\Delta f)(0) - \varepsilon) \pi \rho^2 \leq \int_{B(\rho)} \Delta f \leq ((\Delta f)(0) + \varepsilon) \pi \rho^2$$

for $\rho$ sufficiently small. It follows that

$$(\Delta f)(0) - \varepsilon \leq \frac{4}{r^2} \left(\frac{1}{2\pi} \int_{\theta=0}^{2\pi} f(re^{\sqrt{-1}\theta}) d\theta - f(0)\right) \leq (\Delta f)(0) + \varepsilon$$

for $r$ sufficiently small. □

*Notation.* For tangent vectors $\sigma$ and $\tau$ of a Kähler manifold $X$ of complex dimension $n$ the bisectional curvature for $\sigma$ and $\tau$ is

$$R(\sigma, \tau, \tau, \sigma) + R(\sigma, J\tau, J\tau, \sigma),$$

where $R(\cdot, \cdot, \cdot, \cdot)$ is the Riemann curvature tensor and $J$ is the complex structure operator of the tangent bundle $T_X$ of $X$. When $\sigma = 2\operatorname{Re}\xi$ and $\tau = 2\operatorname{Re}\eta$ for tangent vectors $\xi = \xi^\alpha \frac{\partial}{\partial z^\alpha}$ and $\eta = \eta^\alpha \frac{\partial}{\partial z^\alpha}$, the bisectional curvature for $\sigma$ and $\tau$ is equal to

$$4 R_{\alpha\bar{\beta}\gamma\bar{\delta}} \xi^\alpha \overline{\xi^\eta} \eta^\gamma \overline{\eta^\delta},$$

where $R_{\alpha\bar{\beta}\gamma\bar{\delta}}$ are the components of the Riemann curvature tensor with respect to the local holomorphic coordinates $z^1, \cdots, z^n$ of $X$. When $X$ is compact, for a closed subset $A$ of $X$ we denote by $\operatorname{dist}_A$ the function on $X$ whose value at a point $P$ is the distance from $P$ to $A$.

PROPOSITION 2.2. *Let $X$ be a compact Kähler manifold and $M$ be a $C^m$ Levi-flat real hypersurface in $X$ with $m \geq 3$. Let $P_0$ be a point of $X - M$ such that the distance function $\operatorname{dist}_M$ to $M$ is $C^{m-1}$ on some open neighborhood of $P_0$ and such that the shortest geodesic from $P_0$ to $M$ is represented by a smooth curve $\Phi : [0, \ell] \to X$ parametrized by arc-length with $\Phi(0) = P_0$ and $\Phi(\ell) \in M$. Let $z^j = x^j + \sqrt{-1} y^j$ $(1 \leq j \leq n)$ be a local holomorphic coordinate system at $P_0$ with $\frac{\partial}{\partial x^j}, \frac{\partial}{\partial y^j}$ perpendicular to $\Phi$ at $P_0$ for $1 \leq j < n$. Let $T = (d\Phi)\left(\frac{\partial}{\partial t}\right)$ be the unit tangent vector of the geodesic $\Phi$ (where $t$ is the coordinate of $(0, \ell]$)* and let $\sigma_j$ and $\tau_j = J\sigma_j$ be respectively the parallel vector fields along $\Phi$ whose values at $P_0$ are $\frac{\partial}{\partial x^j}$ and $\frac{\partial}{\partial y^j}$ ($1 \leq j < n$). Then for $1 \leq j < n$ at $P_0$

$$4\frac{\partial^2}{\partial z^j \partial \overline{z^j}} \mathrm{dist}_M \leq -\int_{t=0}^{\ell} (R(\sigma_j, T, T, \sigma_j) + R(\tau_j, T, T, \tau_j))\, dt.$$

In particular, if there exists a positive number $\kappa$ such that for any mutually orthogonal unit tangent vectors $\xi_1, \cdots, \xi_q$ of $X$ of type $(1,0)$ and any unit tangent vector $\eta$ of $X$ of type $(1,0)$ one has

$$\sum_{j=1}^{q} R_{\alpha\bar{\beta}\gamma\bar{\delta}} \xi_j^\alpha \overline{\xi_j^\eta} \eta^\gamma \overline{\eta_j^\delta} \geq \kappa,$$

then the sum of any $q$ eigenvalues of $\sqrt{-1}\partial\bar{\partial}(-\log \mathrm{dist}_M)$ with respect to the Kähler metric is no less than $\kappa$ at points of $X$ near $M$.

*Proof.* Fix $1 \leq j < n$. We construct a $C^{m-1}$ map

$$\tilde{\Phi} : [0, \ell] \times (-\varepsilon, \varepsilon) \times (-\varepsilon, \varepsilon) \to X$$

$$(t, u, v) \mapsto \gamma(t, u, v)$$

so that $\tilde{\Phi}(t, 0, 0) = \Phi(t)$ and the differential $d\tilde{\Phi}$ of $\tilde{\Phi}$ maps $\frac{\partial}{\partial u}$ and $\frac{\partial}{\partial v}$ to $\sigma_j$ and $\tau_j$ respectively. Let $L(u, v)$ be the arc-length of the curve

$$\begin{array}{rcl}[0, \ell] & \to & X \\ t & \mapsto & \tilde{\Phi}(t, u, v).\end{array}$$

Let $g(\cdot, \cdot)$ be the Kähler metric of $X$ and let $Q_0 = \Phi(\ell)$. By the second variation formula

$$\left(\frac{\partial^2}{\partial u^2} + \frac{\partial^2}{\partial v^2}\right) L(u, v)\Big|_{(u,v)=(0,0)} = \left[g(\nabla_{\sigma_j}\sigma_j, T) + g(\nabla_{\tau_j}\tau_j, T)\right]_{P_0}^{Q_0}$$
$$- \int_{t=0}^{\ell} (R(\sigma_j, T, T, \sigma_j) + R(\tau_j, T, T, \tau_j))\, dt.$$

(see e.g., [Fra, formula (7), p. 171]). Since both $\{z_j = \mathrm{constant}\}$ and $M$ are Levi-flat, it follows that

$$g(\nabla_{\sigma_j}\sigma_j, T) + g(\nabla_{\tau_j}\tau_j, T)$$

vanishes both at $P_0$ and at $Q_0$ (see e.g., [Fra, formula (8), p. 171]). Hence

$$\left(\frac{\partial^2}{\partial u^2} + \frac{\partial^2}{\partial v^2}\right) L(u, v)\Big|_{(u,v)=(0,0)} = -\int_{t=0}^{\ell} (R(\sigma_j, T, T, \sigma_j) + R(\tau_j, T, T, \tau_j))\, dt.$$

Since

$$\begin{array}{rcl}L(u, v) & \geq & \mathrm{dist}_M(\tilde{\Phi}(0, u, v)) \text{ for all } (u, v), \\ L(0, 0) & = & \mathrm{dist}_M(\tilde{\Phi}(0, 0, 0)),\end{array}$$

it follows from Lemma (2.1) that

$$\frac{\partial^2}{\partial z^j \partial \overline{z^j}} \operatorname{dist}_M \leq \frac{\partial^2}{\partial z^j \partial \overline{z^j}} L$$

at $P_0$ and

$$4 \frac{\partial^2}{\partial z^j \partial \overline{z^j}} \operatorname{dist}_M \leq -\int_{t=0}^{\ell} \left( R(\sigma_j, T, T, \sigma_j) + R(\tau_j, T, T, \tau_j) \right) dt$$

at $P_0$. $\square$

Similar second-variation arguments in somewhat different settings were given in [T], [E], [Suz].

We now recall the following definition which was introduced in [Siu2, p. 88].

*Definition* 2.3. The bisectional curvature of a Kähler manifold $X$ is said to be *strongly s-nondegenerate* when the following holds. If $k$ and $\ell$ are positive integers, and $\xi_{(1)}, \cdots, \xi_{(k)}$ (respectively $\eta_{(1)}, \cdots, \eta_{(\ell)}$) are **C**-linearly independent tangent vectors of $X$ of type $(1,0)$ such that

$$R_{\alpha\bar{\beta}\gamma\bar{\delta}} \xi_{(\mu)}^{\alpha} \overline{\xi_{(\mu)}^{\beta}} \eta_{(\nu)}^{\gamma} \overline{\eta_{(\nu)}^{\delta}} = 0$$

for $1 \leq \mu \leq k$ and $1 \leq \nu \leq \ell$, then $k + \ell \leq s$. When the condition is satisfied only for the special case of $k = 1$, we say that the bisectional curvature of $X$ is $s$-nondegenerate. The smallest $s$ so that the bisectional curvature of $X$ is (strongly) $s$-nondegenerate is called the *degree of the (strong) nondegeneracy* of the bisectional curvature of $X$.

Clearly, the degree of strong nondegeneracy of the bisectional curvature is no smaller than the degree of nondegeneracy of the bisectional curvature. However, from the computations carried out for irreducible compact Hermitian symmetric manifolds [C-V], [B], [Siu1],[Z], it turns out that the degree of strong nondegeneracy of the bisectional curvature of an irreducible compact Hermitian symmetric manifold is always equal to the degree of nondegeneracy of its bisectional curvature, which is given as follows.

(1) The degree of (strong) nondegeneracy of the bisectional curvature of $\mathrm{U}(m+n)/\mathrm{U}(m) \times \mathrm{U}(n)$ is $(m-1)(n-1) + 1$.

(2) The degree of (strong) nondegeneracy of the bisectional curvature of $\mathrm{SO}(2n)/\mathrm{U}(n) \times \mathrm{U}(n)$ is $\frac{1}{2}(n-2)(n-3) + 1$.

(3) The degree of (strong) nondegeneracy of the bisectional curvature of $\mathrm{Sp}(n)/\mathrm{U}(n)$ is $\frac{1}{2}n(n-1) + 1$.

(4) The degree of (strong) nondegeneracy of the bisectional curvature of $\mathrm{SO}(m+2)/\mathrm{SO}(m) \times \mathrm{SO}(2)$ is 2.

(5) The degree of (strong) nondegeneracy of the bisectional curvature of $E_6/\mathrm{Spin}(10) \times \mathrm{SO}(2)$ is 6.

(6) The degree of (strong) nondegeneracy of the bisectional curvature of $E_7/E_6 \times \mathrm{SO}(2)$ is 11.

PROPOSITION 2.4. *Let $X$ be an irreducible compact Hermitian symmetric manifold of complex dimension $n$ and let $s$ be the degree of the nondegeneracy of the bisectional curvature of $X$. Let $M$ be a $C^m$ Levi-flat real hypersurface in $X$ with $m \geq 3$. Let $\omega_0$ be the standard Kähler form on $X$. Let $K$ be a compact subset of $X - M$ so that from any point of $X - (M \cup K)$ to $M$ there exists a unique minimal geodesic (and in the last statement of Proposition 2.2 points of $X$ near $M$ mean points in $X - (M \cup K)$). Then the following conclusions hold.*

(1) *$X - M$ is Stein.*

(2) *$-\mathrm{dist}_M$ is $C^{m-1}$ and $-\log \mathrm{dist}_M$ is weakly plurisubharmonic on $X - (M \cup K)$.*

(3) *There exists a positive number $\kappa_0$ with the property that the sum of at least $s$ eigenvalues of $\sqrt{-1}\partial\bar\partial(-\log \mathrm{dist}_M)$ is no less than $\kappa_0$ at every point of $X - (M \cup K)$.*

(4) *Let $A$ be a positive number less than the distance between $K$ and $M$. Let $\chi$ be a $C^\infty$ function on the real line $\mathbf{R}$ with $\chi' \geq 0$ and $\chi'' \geq 0$ everywhere such that $\chi \equiv 0$ on $(-\infty, -\log A]$ and $\chi(\lambda) = \lambda$ for $\lambda \geq -\log \frac{A}{2}$. Let $\psi$ be a $C^\infty$ function on $X - M$ which is strictly plurisubharmonic on $\{\mathrm{dist}_M \geq \frac{A}{3}\}$ and whose support is contained in $\{\mathrm{dist}_M \geq \frac{A}{4}\}$. Then the Kähler metric defined by*

$$\omega := \omega_0 + \sqrt{-1}\partial\bar\partial\left(\chi \circ (-\log \mathrm{dist}_M)\right)$$

*is a complete Kähler metric on $X - M$. Moreover, for a sufficiently small positive number $\varepsilon$ there exists a positive number $\kappa$ such that the sum of at least $s$ eigenvalues of*

$$\sqrt{-1}\partial\bar\partial\left(\varepsilon\psi + \chi \circ (-\log \mathrm{dist}_M)\right)$$

*with respect to the Kähler form $\omega$ is no less than $\kappa$.*

(5) *Let $1 \leq q < n - s$ and $\ell \geq 1$ be integers. Let $v$ be a measurable $\bar\partial$-closed $(0, q)$-form on $X - M$ such that $\frac{1}{\mathrm{dist}_M^{\ell+2}} v$ is $L^2$ on $X$ with respect to $\omega_0$. Then there exists a measurable $(0, q-1)$-form $u$ on $X - M$ such that $\bar\partial u = v$ on $X - M$ and $\frac{1}{\mathrm{dist}_M^\ell} u$ is $L^2$ on $X$ with respect to $\omega_0$. Moreover, the solution $u$ is $C^{k+1}$ on any open subset of $X - M$ where $v$ is $C^k$ if $k + 1 \leq m - 1$.*

*Proof.* Conclusion (1) follows from [H]. Conclusion (2) follows from Proposition 2.2 (or from techniques of [E], [H], [Suz], [T]). Conclusion (3) follows from Proposition 2.2 and Definition 2.3. To prove Conclusion (4) we first observe that, since both $\chi'$ and $\chi''$ are everywhere nonnegative and $\chi \equiv 0$ on $(-\infty, -\log A]$, it follows from Conclusion (2) that $\chi \circ (-\log \mathrm{dist}_M)$ is weakly plurisubharmonic everywhere on $X - M$. Hence on $X - M$ the $(1,1)$-form

$$\omega = \omega_0 + \sqrt{-1}\left(\chi \circ (-\log \mathrm{dist}_M)\right)$$

defines a Kähler metric. The completeness of $X - M$ with respect to the Kähler metric $\omega$ results from considering the growth behavior at the origin of the Laplacian of the logarithm of the absolute value of the real part of a complex variable.

To check the remaining statement in Conclusion (4) about the lower bound of a sum of $s$ eigenvalues, first consider a point $P_0$ in $\{0 < \mathrm{dist}_M \leq \frac{A}{4}\}$. Let $\lambda_1, \cdots, \lambda_n$ be the eigenvalues of $\sqrt{-1}\partial\bar{\partial}\left(\chi \circ (-\log \mathrm{dist}_M)\right)$ with respect to $\omega_0$. It follows from Conclusions (2) and (3) that each $\lambda_j$ is nonnegative and

$$\lambda_{j_1} + \cdots + \lambda_{j_s} \geq \kappa_0$$

for any $1 \leq j_1 < \cdots < j_s \leq n$, which implies that one of

$$\lambda_{j_1}, \cdots, \lambda_{j_s}$$

is at least $\frac{\kappa_0}{s}$. At $P_0$ the eigenvalues of $\sqrt{-1}\partial\bar{\partial}\chi \circ (-\log \mathrm{dist}_M)$ with respect to $\omega$ are $\frac{\lambda_1}{1+\lambda_1}, \cdots \frac{\lambda_n}{1+\lambda_n}$. Since $\frac{a}{1+a} \geq \min\left(\frac{a}{2}, \frac{1}{2}\right)$ for any nonnegative number $a$, it follows that

$$\frac{\lambda_{j_1}}{1+\lambda_{j_1}} + \cdots + \frac{\lambda_{j_s}}{1+\lambda_{j_s}} \geq \min\left(\frac{\kappa_0}{2s}, \frac{1}{2}\right)$$

for any $1 \leq j_1 < \cdots < j_s \leq n$. Next we observe that since $\psi$ is strictly plurisubharmonic on $\{\mathrm{dist}_M \geq \frac{A}{3}\}$, when we consider only points in the compact set $\{\mathrm{dist}_M \geq \frac{A}{3}\}$, Conclusion (4) clearly holds (with $\kappa$ set to be some $\kappa' > 0$). To get Conclusion (4) for points in the remaining compact set $\{\frac{A}{4} \leq \mathrm{dist}_M \leq \frac{A}{3}\}$ we need only choose $\varepsilon$ such that

$$\varepsilon\sqrt{-1}\partial\bar{\partial}\psi \geq \frac{-1}{2s}\min\left(\frac{\kappa_0}{2s}, \frac{1}{2}\right)\omega_0$$

on $\{\frac{A}{4} \leq \mathrm{dist}_M \leq \frac{A}{3}\}$ and then choose $\kappa$ satisfying $0 < \kappa \leq \frac{1}{2}\min\left(\frac{\kappa_0}{2s}, \frac{1}{2s}\right)$ and $\kappa \leq \kappa'$.

To prove Conclusion (5), from the definition of $\omega$ we observe that on $X - M$ we have the following inequalities comparing $\omega$ and $\omega_0$ and their volume forms:

$$\frac{c_1}{\mathrm{dist}_M^2} \text{ volume form of } \omega_0 \leq \text{ volume form of } \omega \leq \frac{c_1}{\mathrm{dist}_M^2} \text{ volume form of } \omega_0$$

for some positive constant $c_1, c_2$, and

$$c_1' \omega_0 \leq \omega \leq \frac{c_2'}{\text{dist}_M^2} \omega_0$$

for some positive constants $c_1', c_2'$. Hence $\frac{1}{\text{dist}_M^\ell} v$ is $L^2$ on $X - M$ with respect to $\omega$. For the trivial line bundle on $X - M$ we introduce the Hermitian metric $e^{-\ell\varphi}$ with

$$\varphi = \varepsilon\psi + \chi \circ (-\log \text{dist}_M).$$

We now use Conclusion (4) and apply Theorem (1.1) to the complete Kähler manifold $X - M$ with Kähler metric $\omega$ and to the trivial line bundle with Hermitian metric $e^{-\ell\varphi}$. The $\bar{\partial}$-closed $(0,q)$-form $v$ on $X - M$ is $L^2$ with respect to the weight $e^{-\ell\varphi}$ and the metric $\omega$. By Theorem (1.1) there exists a measurable $(0, q-1)$-form $u$ on $X - M$ such that $\bar{\partial}u = v$ and $\frac{1}{\text{dist}_M^\ell} u$ is $L^2$ with respect to the weight $e^{-\ell\varphi}$ and the metric $\omega$, which implies that $\frac{1}{\text{dist}_M^\ell} u$ is $L^2$ on $X$ with respect to the metric $\omega_0$ and without any other weight, because of the inequalities comparing $\omega$ and $\omega_0$ and their volume forms. It also follows from Theorem (1.1) that $u$ is $C^{k+1}$ on any open subset of $X - M$ where $v$ is $C^k$ if $k + 1 \leq m - 1$. □

## 3. Integral formula for solving the $\bar{\partial}$-equation for $(0,1)$-forms

We recall here the explicit integral formula for solving the $\bar{\partial}$-equation for $(0,1)$-forms on the ball from the theory of Henkin and Grauert-Lieb (see e.g., [R, Chap. V]). From

$$|\zeta|^2 - |z|^2 = 2\text{Re}\left(\sum_{j=1}^n (\zeta_j - z_j)\bar{\zeta}_j\right) - |\zeta - z|^2$$

it follows that for the ball $B_R$ of radius $R$ in $\mathbf{C}^n$ centered at the origin one has

$$\sum_{j=1}^n (\zeta_j - z_j)\bar{\zeta}_j \neq 0$$

for $\zeta \in \partial B_R$ and $z \in B_R$. Let

$$g_j(\lambda, \zeta, z) = \lambda \frac{\bar{\zeta}_j}{\sum_{k=1}^n (\zeta_k - z_k)\bar{\zeta}_k} + (1-\lambda)\frac{\bar{\zeta}_j - \bar{z}_j}{|\zeta - z|^2}.$$

From

$$\sum_{j=1}^n (\zeta_j - z_j)\, g_j(\lambda, \zeta, z) = 1$$

it follows that
$$\sum_{j=1}^{n} (\zeta_j - z_j)(\bar{\partial}_\zeta + d_\lambda) g_j(\lambda, \zeta, z) = 0$$

and
$$\bigwedge_{j=1}^{n} ((\bar{\partial}_\zeta + d_\lambda) g_j(\lambda, \zeta, z)) = 0.$$

Let
$$c_n = \frac{1}{n} \left( \frac{\sqrt{-1}}{2\pi} \right) (-1)^{\frac{n(n+1)}{2}}$$

and
$$\Omega(\lambda, \zeta, z) = c_n \left( \sum_{j=1}^{n} (-1)^{j-1} g_j \bigwedge_{1 \le k \le n, k \ne j} ((\bar{\partial}_\zeta + d_\lambda) g_k(\lambda, \zeta, z)) \right) \bigwedge_{k=1}^{n} (d\zeta_k).$$

Then
$$(d_\zeta + d_\lambda) \Omega(\lambda, \zeta, z) = 0$$

for $\zeta \ne z$. The restriction of $\Omega$ to $\{\lambda = 0\}$ is equal to the Bochner-Martinelli kernel
$$K_{BM}(\zeta, z) = c_n \left( \sum_{j=1}^{n} (-1)^{j-1} \frac{\bar{\zeta}_j - \bar{z}_j}{|\zeta - z|^2} \bigwedge_{1 \le k \le n, k \ne j} \bar{\partial}_\zeta \frac{\bar{\zeta}_j - \bar{z}_j}{|\zeta - z|^2} \right) \bigwedge_{k=1}^{n} (d\zeta_k).$$

From
$$\left( \frac{\sqrt{-1}}{2\pi} \partial_\zeta \bar{\partial}_\zeta \log |\zeta - z|^2 \right)^n = \delta_z$$

it follows that on $\{\lambda = 0\}$ one has
$$(d_\zeta + d_\lambda) \Omega(\lambda, \zeta, z) = \delta_z,$$

where $\delta_z$ is the Kronecker delta at $z$.

For a function $f(\zeta)$ we apply Stokes' theorem to
$$(d_\zeta + \lambda)(\Omega(\lambda, \zeta, z) f(\zeta))$$

on the bordered manifold defined by
$$\{(\lambda, \zeta) \mid \lambda = 0, \zeta \in \overline{B_R}\} \cup \{(\lambda, \zeta) \mid 0 \le \lambda \le 1, \zeta \in \partial B_R\}$$

and get

$$f(z) = \int_{\zeta \in \partial B_R} \Omega(\lambda, \zeta, z)|_{\lambda=1} f(\zeta) + \int_{\zeta \in B_R} K_{BM}(\zeta, z) \wedge \bar{\partial} f(\zeta)$$
$$+ \int_{\zeta \in \partial B_R} \left( \int_{\lambda=0}^{1} \Omega(\lambda, \zeta, z) \right) \wedge \bar{\partial} f(\zeta).$$

Since $\Omega(\lambda, \zeta, z)|_{\lambda=1}$ is holomorphic in $z \in B_R$ it follows that for a $\bar{\partial}$-closed $(0,1)$-form $\alpha$ on $\overline{B_R}$ the function

$$v_R(z) := \int_{\zeta \in B_R} K_{BM}(\zeta, z) \wedge \alpha + \int_{\zeta \in \partial B_R} \left( \int_{\lambda=0}^{1} \Omega(\lambda, \zeta, z) \right) \wedge \alpha$$

satisfies $\bar{\partial} v_R = \alpha$. Let $\rho(\lambda)$ be a nonnegative $C^\infty$ function supported on $(r_1^2, r_2^2)$ with $0 < r_1 < r_2 < R$ such that $\int_{\lambda=r_1^2}^{r_2^2} \rho(\lambda) d\lambda = 1$. Then

$$\tilde{v}(z) := \int_{r=r_1}^{r_2} \rho(r^2) v_r(z) d(r^2)$$
$$= \int_{|\zeta|<r_1} K_{BM}(\zeta, z) \wedge \alpha + \int_{r=r_1}^{r_2} \left( \int_{r_1<|\zeta|<r} K_{BM}(\zeta, z) \wedge \alpha \right) \rho(r^2) d(r^2)$$
$$+ \int_{r=r_1}^{r_2} \left( \int_{\zeta \in \partial B_r} \left( \int_{\lambda=0}^{1} \Omega(\lambda, \zeta, z) \right) \wedge \alpha \right) \rho(r^2) d(r^2)$$

satisfies $\bar{\partial} \tilde{v}(z) = \alpha(z)$ for $|z| < r_1$. We rewrite $\tilde{v}(z)$ in the form

$$\tilde{v}(z) = \int_{|\zeta|<r_1} K_{BM}(\zeta, z) \wedge \alpha + \int_{r_1<|\zeta|<R} \tilde{K}(\zeta, z) \wedge \alpha,$$

where $\tilde{K}(\zeta, z)$ is $C^\infty$ for $|z| < r_1 < |\zeta| < R$. The Bochner-Martinelli kernel $K_{BM}(\zeta, z)$ is of the form

$$K_{BM}(\zeta, z) = \frac{K_1(\zeta, z)}{|\zeta - z|^{4n-2}},$$

where $K_1(\zeta, z)$ is $C^\infty$ and satisfies $|K_1(\zeta, z)| \leq C|\zeta - z|^{2n-1}$ for some constant $C$. For our purpose we need the bound and the convergence behavior of the integral only for $z \in B_{r_1}$. Only the contribution from the Bochner-Martinelli kernel $K_{BM}(\zeta, z)$ needs to be considered.

## 4. Regularity of $\bar{\partial}_b$

LEMMA 4.1. *Let $p$ and $\ell$ be positive integers such that $\ell > \frac{3n}{2} - 1$ and $p < \ell - \frac{3n}{2} + 1$. Let $W$ be an open neighborhood of the origin in $\mathbf{C}^n$ and $f$ be a $C^p$ real-valued function on an open neighborhood $U$ of the closure of $W$ in $\mathbf{C}^n$ such that $f(0) = 0$ and $df$ is nowhere zero on $U$. Let $S_\sigma$ be the real submanifold of $W$ defined by $f = \sigma$. Assume that for some $\sigma^* > 0$ there is a $C^p$ family of $C^p$ diffeomorphisms $\varphi_\sigma$ parametrized by $\sigma \in (-\sigma^*, \sigma^*)$ such that for each $\sigma$ the diffeomorphism $\varphi_\sigma$ maps $S_0$ diffeomorphically onto $S_\sigma$ and $\varphi_0$ is the identity.*

(1) Let $\alpha > 1$ and $u$ be a measurable function on $W$ such that $\frac{u}{|f|^\alpha}$ is $L^1$. Then there exists a subsequence $\sigma_\nu \to 0$ such that $\int_{S_{\sigma_\nu}} |u| \to 0$ as $\nu \to \infty$.

(2) Let $u(z)$ be $L^2$ on $W$ and let $K(z,\zeta)$ be a function on $U$ of the form $\frac{K_1(z,\zeta)}{|z-\zeta|^{4n-2}}$ such that $K_1(z,\zeta)$ is $C^p$ and $|K_1(z,\zeta)| \le C|z-\zeta|^{2n-1}$ on $U$ for some constant $C$. Let $r > 0$ such that the ball $B(r)$ of radius $r$ centered at the origin is contained in $W$ and let
$$v(z) = \int_{\zeta \in B(r)} K(z,\zeta) u(\zeta) d\lambda(\zeta),$$
where $d\lambda(\zeta)$ is the Euclidean volume form in the variable $\zeta$. Then each $v|_{S_\sigma}$ is $L^1$ and approaches $v|_{S_0}$ in the $L^1$ norm as $\sigma \to 0$ (in the sense that $v \circ \varphi_\sigma$ approaches $v$ on $S_0$ in the $L^1$ norm on $S_0$ as $\sigma \to 0$).

(3) If $u$ and $v$ are as in (2) with the additional assumption that $\frac{u}{|f|^\ell}$ is $L^2$ on $U$, then the $L^1$ function $v|_{S_0}$ on $S_0$ is $C^p$.

(4) Let $\omega$ be a $\bar\partial$-closed $(0,1)$-current on $U$ which is $C^{p-1}$ on $U - S_0$ such that $\frac{\omega}{|f|^\ell}$ is $L^2$ on $U$. Let $v$ be a generalized function on $W$ with $\bar\partial v = \omega$ on $W$. Then $v|_{S_\sigma}$ approaches some $C^p$ function $v_0$ on $S_0$ in the $L^1$ norm as $\sigma \to 0$ and $\bar\partial v_0 = 0$ on $S_0$. Here $v|_{S_\sigma}$ approaching $v_0$ in the $L^1$ norm means that $v \circ \varphi_\sigma$ approaches $v_0$ in the $L^1$ norm on $S_0$ as $\sigma \to 0$.

*Proof.* At an arbitrary point $P_0$ of $S_0$ we take a $C^p$ local coordinate system $g_1, \cdots, g_n$ at a point of $S_0$ with $\operatorname{Im} g_n = f$ such that, for some $\sigma_{P_0} > 0$ and some open neighborhood $G_{P_0}$ of $P_0$ in $S_0$, the diffeomorphism $\varphi_\sigma$ maps the point in $G_{P_0}$ with coordinates $(g_1, \cdots, g_{n-1}, \operatorname{Re} g_n)$ to the point in $S_\sigma$ with coordinates $(g_1, \cdots, g_{n-1}, \operatorname{Re} g_n + \sqrt{-1}\sigma)$ for $|\sigma| < \sigma_{P_0}$ with respect to the coordinate system $(g_1, \cdots, g_n)$. By replacing our coordinate system $(z_1, \cdots, z_n)$ by $(g_1, \cdots, g_n)$ and replacing $W$ by a smaller neighborhood of $0$ we can assume without loss of generality, for the proofs of (1), (2), and (3), that $f$ equals the imaginary part $y_n$ of $z_n$ and that the diffeomorphism $\varphi_\sigma$ maps the point $(z_1, \cdots, z_{n-1}, \operatorname{Re} z_n)$ to the point $(z_1, \cdots, z_{n-1}, \operatorname{Re} z_n + \sqrt{-1}\sigma)$ in the coordinate system $(z_1, \cdots, z_n)$.

(1) Let $U(\sigma) = \int_{\{y_n = \sigma\}} \left|\frac{u}{y_n^\alpha}\right|$. There exists no $\varepsilon > 0$ such that $U(\sigma) \ge \frac{1}{\sigma}$ for all $0 < \sigma < \varepsilon$; otherwise we have the contradictory conclusion
$$\infty \le \int_{\sigma=0}^\varepsilon \frac{d\sigma}{\sigma} \le \int_{\sigma=0}^\varepsilon U(\sigma) d\sigma = \int_W \left|\frac{u}{y_n^\alpha}\right| < \infty.$$
Thus there exists $0 < \sigma_\nu < \frac{1}{\nu}$ such that $U(\sigma_\nu) < \frac{1}{\sigma_\nu}$ for all $\nu$. We have
$$\int_{\{y_n = \sigma_\nu\}} |u| = \sigma_\nu^\alpha U(\sigma_\nu) \le \sigma_\nu^{\alpha-1} \to 0$$
as $\nu \to \infty$.

(2) Write $z = (z_1, \cdots, z_n)$ and $\zeta = (\zeta_1, \cdots, \zeta_n)$. Let $z_\nu = x_\nu + \sqrt{-1} y_\nu$ and $\zeta_\nu = \xi_\nu + \sqrt{-1} \eta_\nu$ $(1 \leq \nu \leq n)$. Let

$$v_\delta(z) = \int_{\zeta \in B(r) \cap \{|\eta_n| < \delta\}} K(z, \zeta) u(\zeta) d\lambda(\zeta),$$

$$w_\delta(z) = \int_{\zeta \in B(r) \cap \{|\eta_n| \geq \delta\}} K(z, \zeta) u(\zeta) d\lambda(\zeta),$$

so that $v(z) = v_\delta(z) + w_\delta(z)$.

First we claim that for any $\varepsilon > 0$ there exists some $\delta_0 > 0$ such that

$$\int_{\{y_n = \sigma\}} |v_\delta(z)| < \varepsilon$$

for $0 < \delta \leq \delta_0$ and for all $\sigma$. We use $C_1, C_2, C_3$ to denote positive constants. Let

$$\tau = \left( |x_n - \xi_n|^2 + \sum_{\nu=1}^{n-1} |z_\nu - \zeta_\nu|^2 \right)^{\frac{1}{2}}.$$

Let $R = 2r$. Using the fact that the Jacobian in the polar coordinate with radius $\tau$ is $\tau^{2n-2}$, we have

$$\int_{\{\tau \leq R, y_n = \sigma\}} \frac{1}{|z - \zeta|^{2n-1}} dx_1 \cdots dx_n dy_1 \cdots dy_{n-1}$$

$$= \int_{\{\tau \leq R, y_n = \sigma\}} \frac{1}{(\tau + |y_n - \eta_n|)^{2n-1}} dx_1 \cdots dx_n dy_1 \cdots dy_{n-1}$$

$$\leq C_1 \int_{\tau=0}^{R} \frac{1}{(\tau + |y_n - \eta_n|)^{2n-1}} \tau^{2n-2} d\tau$$

$$\leq C_1 \int_{\tau=0}^{R} \frac{d\tau}{\tau + |y_n - \eta_n|}$$

$$= C_1 \log\left(1 + \frac{R}{|y_n - \eta_n|}\right) \leq C_2 \log \frac{1}{|y_n - \eta_n|}.$$

Thus

$$\int_{\{y_n = \sigma\}} |v_\delta(z)| dx_1 \cdots dx_n dy_1 \cdots dy_{n-1}$$

$$\leq C C_2 \int_{\{|\zeta| < r, |\eta_n| < \delta\}} \left| \log \frac{1}{|\sigma - \eta_n|} \right| |u(\zeta)| d\xi_1 \cdots d\xi_n d\eta_1 \cdots d\eta_n$$

$$\leq C C_2 \left( \int_{\{|\zeta| < r, |\eta_n| < \delta\}} \left| \log \frac{1}{|\sigma - \eta_n|} \right|^2 d\xi_1 \cdots d\xi_n d\eta_1 \cdots d\eta_n \right)^{\frac{1}{2}}$$

$$\cdot \left( \int_{\{|\zeta| < r, |\eta_n| \leq \delta\}} |u(\zeta)|^2 d\xi_1 \cdots d\xi_n d\eta_1 \cdots d\eta_n \right)^{\frac{1}{2}}$$

$$\leq C C_3 \left( \int_{\{|\zeta| < r, |\eta_n| \leq \delta\}} |u(\zeta)|^2 d\xi_1 \cdots d\xi_n d\eta_1 \cdots d\eta_n \right)^{\frac{1}{2}},$$

which implies the claim, because
$$\int_{\{|\zeta|<r,|\eta_n|\leq\delta\}} |u(\zeta)|^2 \to 0$$
as $\delta \to 0$.

Now for fixed $\delta = \delta_0$ the function $w_\delta(z)$ is $C^p$ in $z$ for $|y_n| < \delta$. Hence there exists $\sigma_0 > 0$ such that the $L^1$-norm of $w_\delta|_{y_n=\sigma} - w_\delta|_{y_n=0}$ as a function of the variables $z_1, \cdots, z_{n-1}, x_n$ is less than $\varepsilon$ for $|\sigma| < \sigma_0$. Thus the $L^1$-norm of $v|_{y_n=\sigma} - v|_{y_n=0}$ as a function of the variables $z_1, \cdots, z_{n-1}, x_n$ is less than $3\varepsilon$ for $|\sigma| < \sigma_0$, which proves Conclusion (2).

(3) To verify that $v|_{y_n=0}$ is $C^p$, we write
$$v(z) = \lim_{\delta\to 0} \int_{z\in B(r)\cap\{|\eta_n|\geq\delta\}} K(z,\zeta)u(\zeta)d\lambda(\zeta).$$

On the manifold $\{y_n = 0\}$, to justify the commutation of differentiation along a real direction and integration in the above improper integral, one need only check that the result of the differentiation is absolutely convergent. Any derivative $D^p K(z,\zeta)$ of the kernel $K(z,\zeta)$ to order $p$ (along the direction in $\{y_n = 0\}$) satisfies
$$|D^p K(z,\zeta)| \leq \frac{C_4}{|z-\zeta|^{2n+p-1}}$$
for some constant $C_4$. Since $\frac{u(z)}{|y_n|^\ell}$ is $L^2$, it follows that
$$\int_{\zeta\in B_r} |D^p K(z,\zeta) u(\zeta)| d\lambda(\zeta)$$
$$\leq C_4 \left(\int_{\zeta\in B_r} \frac{|\eta_n|^{2\ell}}{|z-\zeta|^{2(2n+p-1)}} d\lambda(\zeta)\right)^{\frac{1}{2}} \left(\int_{\zeta\in B_r} \left|\frac{u(\zeta)}{\eta_n^\ell}\right|^2 d\lambda(\zeta)\right)^{\frac{1}{2}}$$
which is finite for $2\ell > 2(2n+p-1) - n = 3n+2p-2$ or $\ell > \frac{3n}{2} + p - 1$.

(4) To prove Conclusion (4), we use a new $C^p$ local coordinate system of $\mathbf{C}^n$ to straighten out the family $\{S_\sigma\}$ and to express the integral formula in Section 3 for solving the $\bar\partial$-equation and the $\bar\partial_b$ operator on each $S_\sigma$ in terms of the new $C^p$ local coordinate system. We use pointwise $\mathbf{C}$-linearly independent local $C^{p-1}$ $(1,0)$-forms $\theta_1, \cdots, \theta_n$ on $\mathbf{C}^n$ such that the pullback of $\theta_n$ to $T^{1,0}_{S_\sigma}$ is identically zero for all $\sigma$. We write $\omega = \sum_{\nu=1}^n u^{(\nu)} \bar\theta_\nu$ and write the $\bar\partial_b$ operator on $S_\sigma$ as $\sum_{\nu=1}^{n-1} (\bar\theta_\nu|_{S_\sigma}) \otimes \xi_\sigma^{(\nu)}$ for some vector fields $\xi_\sigma^{(\nu)}$ on $S_\sigma$. More precisely, at an arbitrary point $P_0$ of $S_0$, because of the integral formulas in Section 3, we can choose a $C^p$ local coordinate system $g_1, \cdots, g_n$ at a point of $S_0$ adapted to the $C^p$ family of real hypersurfaces $S_\sigma$ in the sense that the following conditions hold.

(i) $\operatorname{Im} g_n = f$.

(ii) The diffeomorphism $\varphi_\sigma$ maps the point $(z_1, \cdots, z_{n-1}, \operatorname{Re} z_n)$ to the point $(z_1, \cdots, z_{n-1}, \operatorname{Re} z_n + \sqrt{-1}\,\sigma)$ in the coordinate system $(z_1, \cdots, z_n)$.

(iii) There exist functions $K^{(\nu)}(z,\zeta)$ $(1 \leq \nu \leq n)$ on $U$ of the form $\frac{K_1^{(\nu)}(z,\zeta)}{|z-\zeta|^{4n-2}}$ such that $K_1^{(\nu)}(z,\zeta)$ is $C^p$ and $\left|K_1^{(\nu)}(z,\zeta)\right| \leq C|z-\zeta|^{2n-1}$ on $U$ for some constant $C$.

(iv) There exists a holomorphic function $F$ on $U$ and there exist $C^p$ functions $u^{(\nu)}$ $(1 \leq \nu \leq n)$ on $U - S_0$ constructed from $\omega$ such that each $\frac{u^{(\nu)}}{|f|^\ell}$ is $L^2$ and
$$v(z) = F(z) + \sum_{\nu=1}^{n} \int_{\zeta \in B(r)} K^{(\nu)}(z,\zeta) u^{(\nu)}(\zeta).$$

(v) There exist $C^{p-1}$ families of $C^{p-1}$ vector fields $\xi_\sigma^{(\nu)}$ on $S_\sigma$ $(1 \leq \nu < n)$ with parameter $\sigma$ such that $\xi_\sigma^{(\nu)}(v|_{S_\sigma}) = u^{(\nu)}$ on $S_\sigma$ for $\sigma \neq 0$, where $\xi_\sigma^{(\nu)}(v|_{S_\sigma})$ means the result obtained by applying the vector field $\xi_\sigma^{(\nu)}$ to $v|_{S_\sigma}$.

(vi) The condition $\bar{\partial}_b(v|_{S_0}) = 0$ on $S_0$ in the sense of distributions is equivalent to $\xi_0^{(\nu)}(v|_{S_0}) = 0$ on $S_0$ $(1 \leq \nu < n)$ in the sense of distributions.

By Conclusion (2), $v|_{S_\sigma}$ approaches $v|_{S_0}$ in the $L^1$ norm as $\sigma \to 0$. Since $\xi_\sigma^{(\nu)}$ is a $C^{p-1}$ family with parameter $\sigma$, it follows that the distribution $\xi_\sigma^{(\nu)}(v|_{S_\sigma})$ approaches the distribution $\xi_0^{(\nu)}(v|_{S_0})$ in the space of distributions as $\sigma \to 0$. By Conlusion (1), it follows from the $L^2$ property of $\frac{u^{(\nu)}}{|f|^\ell}$ that there exists a sequence $\sigma_\mu$ $(1 \leq \mu < \infty)$ approaching $0$ such that $u^{(\nu)}|_{S_{\sigma_\mu}}$ approaches $0$ in the $L^1$ norm as $\mu \to \infty$, which implies that $\xi_0^{(\nu)}(v|_{S_0}) = 0$ on $S_0$ in the sense of distributions on $S_0$ and, as a consequence, $\bar{\partial}_b(v|_{S_0}) = 0$ on $S_0$ in the sense of distributions. It follows from Conclusion (3) and the $L^2$ property of $\frac{u^{(\nu)}}{|f|^\ell}$ that $v|_{S_0}$ is $C^p$ on $S_0$. □

## 5. $\partial_b \bar{\partial}_b$-exactness and flat metrics

The following simple lemma is standard. We put it here just to get the precise statement we want concerning differentiability.

LEMMA 5.1 (holomorphic foliation in Levi-flat real hypersurface). *Let $m \geq 2$ be an integer. Let $M$ be a $C^m$ Levi-flat real hypersurface in a complex manifold $X$ defined by the vanishing of a $C^m$ function $f$ on $X$ whose differential $df$ is nowhere zero at points of $M$. Let $J$ be the complex structure operator of the tangent bundle $T_X$ of $X$. Then $J df = 0$ defines a $C^{m-1}$ foliation of $M$ by local regular complex hypersurfaces of $X$.*

*Proof.* The Levi-flatness of $M$ means that $d(Jdf) = 0$ identically on the zero-set of $Jdf$, because, by using local coordinates $z_1, \cdots, z_n$ and $J(dz_\nu) = \sqrt{-1}\,dz_\nu$, computation in terms of $dz_\nu, d\overline{z_\nu}$ gives $d(Jdf) = -2\sqrt{-1}\,\partial\bar{\partial}f$. By using a local basis of 1-forms on $M$ with $Jdf$ as a member, we conclude that $d(Jdf) = (Jdf) \wedge \theta$ for some 1-form $\theta$ on $M$. By Frobenius's theorem the system of distributions defined by $Jdf = 0$ is integrable. The local integral hypersurfaces in $M$ are complex-analytic, because their tangent spaces defined by $df = Jdf = 0$ are stabilized by $J$. Since $f$ is $C^m$, the foliation is $C^{m-1}$ along the directions transversal to the leaves. □

LEMMA 5.2 (extension of cohomology classes). *Let $X$ be a complex manifold and $M$ be a $C^{m+1}$ Levi flat real hypersurface in $X$. Let $\alpha$ be a $C^m$ $\bar{\partial}_b$-closed $(0,1)$-form on $M$ (i.e., $\alpha$ is a $C^m$ function on $\left(T_M^{0,1}\right)^*$). Then there exists a $C^{m-1}$ $(0,1)$-form $\tilde{\alpha}$ on $X$ which extends $\alpha$ as a class (in the sense that the restriction of $\tilde{\alpha}$ to $\left(T_M^{0,1}\right)^*$ agrees with $\alpha + \bar{\partial}_b H$ for some $C^m$ function $H$ on $M$) such that*

(a) *$\bar{\partial}\tilde{\alpha}$ vanishes to order at least $m-2$ at points of $M$, and*

(b) *at every point $P$ of $M$ there exist an open neighborhood $U$ of $P$ in $X$ and a $C^m$ function $h$ on $U$ such that $\tilde{\alpha} - \bar{\partial}h$ vanishes at points of $M \cap U$ to order at least $m-2$.*

*Proof.* The idea of the proof is

(i) that $\alpha$ as a class can be represented by a Čech 1-cocyle of $C^m$ $\bar{\partial}_b$-closed functions with respect to a covering of $M$, and

(ii) that any $C^m$ $\bar{\partial}_b$-closed function on an open subset $W$ of $M$ can be uniquely extended to an $m$-jet extension (in a finite neighborhood of order $m$ of $W$ in $X$) whose $\bar{\partial}$ vanishes to order $m-1$ at points of $W$.

To give the details of the proof, we cover $M$ by local coordinate charts $W_j$ such that each $W_j$ is Cauchy-Riemann equivalent to $G_j \times I_j$ by a $C^{m+1}$ diffeomorphism, where $G_j$ is an open subset of $\mathbf{C}^{n-1}$ and $I_j$ is an open interval in $\mathbf{R}$. We choose an open subset $\tilde{W}_j$ in $X$ such that $M \cap \tilde{W}_j = W_j$. Then $\alpha|W_j$ can be written as $\bar{\partial}_b h_j$ on $W_j$ for some $C^m$ function $h_j$ on $W_j$. Let $h_{j\ell} = h_\ell - h_j$ on $W_j \cap W_\ell$. Then $h_{j\ell}$ is $\bar{\partial}_b$-closed on $W_j \cap W_\ell$. The Čech 1-cocyle $\{h_{j\ell}\}_{j,\ell}$ for the covering $\{W_j\}_j$ of $M$ represents the class defined by $\alpha$.

There is a unique $m$-jet extension $\tilde{h}_{j\ell}$ of $h_{j\ell}$ (in a finite neighborhood of order $m$ of $W_j \cap W_\ell$ in $X$) so that $\bar{\partial}\tilde{h}_{j\ell}$ vanishes to order $m-1$ at points of $M$ and is skew-symmetric in $j, \ell$. We can assume that the $m$-jet $\tilde{h}_{j\ell}$ is induced by a $C^m$ function on $\tilde{W}_j \cap \tilde{W}_\ell$ skew-symmetric in $j, \ell$, which we again use the same symbol $\tilde{h}_{j\ell}$ to denote. It follows from the uniqueness that the $m$-jet extension

$$\tilde{h}_{j\ell} + \tilde{h}_{\ell p} + \tilde{h}_{pj}$$

vanishes to order $m$ on $W_j \cap W_\ell \cap W_p$. The class defined by $\{\tilde{h}_{j\ell}\}_{j,\ell}$ will give us our extension by a partition of unity.

We now introduce a partition of unity. We take a $C^\infty$ function $\rho_j$ on $\tilde{W}_j$ with compact support so that $\sum_j \rho_j$ is identically 1 on some open neighborhood $D$ of $M$ in $X$. We define $\tilde{h}_j = \sum_\ell \rho_\ell \tilde{h}_{j\ell}$ to get $\tilde{h}_{j\ell} = \tilde{h}_\ell - \tilde{h}_j$ on $\tilde{W}_\ell \cap \tilde{W}_j \cap D$ up to order $m$ along $W_j \cap W_\ell$. Since $\bar{\partial}\tilde{h}_{j\ell}$ vanishes to order $m-1$ at points of $W_j \cap W_\ell$, it follows that $\bar{\partial}\tilde{h}_j$ agrees with $\bar{\partial}\tilde{h}_\ell$ on $W_j \cap W_\ell$ up to order $m-1$. Thus we can define an $(m-1)$-jet extension $\hat{\alpha}$ of $\alpha$ as a class by setting $\hat{\alpha}$ equal to $\bar{\partial}\tilde{h}_j$ on $\tilde{W}_j$. We then extend $\hat{\alpha}$ to a $C^{m-1}$ $(0,1)$-form $\tilde{\alpha}$ on $X$ so that $\bar{\partial}\tilde{\alpha}$ vanishes to order at least $m-2$ at points of $M$. Clearly from the construction of $\tilde{\alpha}$ the $(0,1)$-form $\tilde{\alpha} - \bar{\partial}\tilde{h}_j$ vanishes to order at least $m-2$ at points of $W_j$ and $\tilde{h}_j$ is $C^m$ on $\tilde{W}_j$.

From $h_\ell - h_j = h_{j\ell} = \tilde{h}_{j\ell} = \tilde{h}_\ell - \tilde{h}_j$ on $W_j \cap W_\ell$ it follows that $h_\ell - \tilde{h}_\ell = h_j - \tilde{h}_j$ on $W_j \cap W_\ell$ and we can define a $C^m$ function $H$ on $M$ by setting $H = \tilde{h}_j - h_j$ on $W_j$. Then on $M$ we have $\tilde{\alpha} = \bar{\partial}_b \tilde{h}_j = \bar{\partial}_b h_j + \bar{\partial}_b H = \alpha + \bar{\partial}_b H$. □

PROPOSITION 5.3. *Let $X$ be an irreducible compact Hermitian symmetric manifold of complex dimension $n$ and let $s$ be the degree of the nondegeneracy of the bisectional curvature of $X$ such that $n - s \geq 2$. Let $m \geq 6$ be an integer. Let $M$ be a $C^{m+1}$ Levi-flat real hypersurface in $X$. Let $\tilde{\alpha}$ be a $C^{m-1}$ $(0,1)$-form $\tilde{\alpha}$ on $X$ such that $\bar{\partial}\tilde{\alpha}$ vanishes to order $m-2$ at points of $M$. Then there exists a $C^{m-1}$ $(0,1)$-form $\beta$ on $X - M$ such that $\frac{1}{\text{dist}_M^{m-4}}\beta$ is $L^2$ on $X$ and $\bar{\partial}\beta = \bar{\partial}\tilde{\alpha}$ on $X$.*

*Proof.* Since $\frac{1}{\text{dist}_M^{m-2}}\bar{\partial}\tilde{\alpha}$ is $L^2$ on $X$, by applying Conclusion (5) of Proposition (2.4) to $v = \bar{\partial}\tilde{\alpha}$ and $\ell = m - 4$ and $q = 2$, we conclude that there exists a $C^{m-1}$ $(0,1)$-form $\beta$ on $X - M$ such that $\bar{\partial}\beta = \bar{\partial}\tilde{\alpha}$ on $X - M$ such that $\frac{1}{\text{dist}_M^{m-4}}\beta$ is $L^2$ on $X$.

Take an arbitrary point $P_0$ of $M$. To verify $\bar{\partial}\beta = \bar{\partial}\tilde{\alpha}$ on some open neighborhood $U$ of $P_0$ in $X$ we choose $U$ so that $M \cap U$ is the zero-set of some $C^{m+1}$ function $f$ on $U$ with $df$ nowhere zero at points of $M \cap U$. We have to check that $(\beta, \bar{\partial}^*\varphi) = (\tilde{\alpha}, \bar{\partial}^*\varphi)$ for any $C^\infty$ $(0,1)$-form $\varphi$ on $U$ with compact support. For $\varepsilon > 0$ let $\psi_\varepsilon$ be a $C^\infty$ function on $U$ such that

(i) $\psi_\varepsilon$ is identically 1 on $U \cap \{\text{dist}_M > \varepsilon\}$,

(ii) $\psi_\varepsilon$ is identically zero on $U \cap \{\text{dist}_M < \frac{\varepsilon}{2}\}$, and

(iii) $|\bar{\partial}\psi_\varepsilon| \leq \frac{C}{\varepsilon}$ on $U$ for some positive constant $C$.

Then clearly $(\beta, \bar\partial^*(\varphi\psi_\varepsilon)) = (\tilde\alpha, \bar\partial^*(\varphi\psi_\varepsilon))$ for any $C^\infty$ $(0,1)$-form $\varphi$ on $U$ with compact support. We now use the fact that

$$f^2 \bar\partial \psi_\varepsilon \to 0$$

in $L^2$ on $U$ as $\varepsilon \to 0$ so that

$$f^2 \bar\partial^*(\varphi\psi_\varepsilon) = f^2 \psi_\varepsilon \bar\partial^* \varphi + \left\langle \varphi, f^2 \bar\partial \psi_\varepsilon \right\rangle \to f^2 \bar\partial^* \varphi$$

in $L^2$ on $U$ as $\varepsilon \to 0$. Let $\tilde\beta = f^{-(m-4)}\beta$ on $U$ which is $L^2$ on $U$. Since $m \geq 6$, it follows that

$$(\beta, \bar\partial^*(\varphi\psi_\varepsilon)) = \left(f^{m-6}\tilde\beta, f^2 \bar\partial^*(\varphi\psi_\varepsilon)\right) \to \left(f^{m-6}\tilde\beta, f^2 \bar\partial^* \varphi\right) = (\beta, \bar\partial^* \varphi)$$

as $\varepsilon \to 0$. Clearly

$$(\tilde\alpha, \bar\partial^*(\varphi\psi_\varepsilon)) = (\bar\partial\tilde\alpha, \varphi\psi_\varepsilon) \to (\bar\partial\tilde\alpha, \varphi) = (\tilde\alpha, \bar\partial^* \varphi)$$

as $\varepsilon \to 0$, because $\tilde\alpha$ is $C^{m-1}$ and $m \geq 6$. From $(\beta, \bar\partial^*(\varphi\psi_\varepsilon)) = (\tilde\alpha, \bar\partial^*(\varphi\psi_\varepsilon))$ it follows that $(\beta, \bar\partial^* \varphi) = (\tilde\alpha, \bar\partial^* \varphi)$. Hence $\bar\partial\beta = \bar\partial\tilde\alpha$. □

PROPOSITION 5.4 (regular solution of $\bar\partial_b$). *Let $X$ be an irreducible compact Hermitian symmetric manifold of complex dimension $n$ and let $s$ be the degree of the nondegeneracy of the bisectional curvature of $X$ such that $n - s \geq 2$. Let $m > \frac{3n}{2} + 4$ be an integer. Let $M$ be a $C^{m+1}$ Levi-flat real hypersurface in $X$. Let $\alpha$ be a $C^m$ $\bar\partial_b$-closed $(0,1)$-form on $M$. Let $p$ be a positive integer less than $m - \frac{3n}{2} - 3$. Then there exists a $C^p$ function $v$ on $M$ such that $\bar\partial_b v = \alpha$ on $M$.*

*Proof.* By Lemma 5.2 there exists a $C^m$ function $H$ on $M$ such that $\alpha + \bar\partial_b H$ can be extended to a $C^{m-1}$ $(0,1)$-form $\tilde\alpha$ on $X$ and $\bar\partial\tilde\alpha$ vanishes to order at least $m - 2$ at points of $M$ and at every point $P$ of $M$ there exist an open neighborhood $U$ of $P$ in $X$ and a $C^m$ function $h$ on $U$ with $\tilde\alpha - \bar\partial h$ vanishing at points of $M \cap U$ to order at least $m - 2$. By Proposition 5.3 there exists a $C^{m-1}$ $(0,1)$-form $\beta$ on $X - M$ such that $\bar\partial\beta = \bar\partial\tilde\alpha$ on $X - M$ and $\frac{1}{\text{dist}_M^{m-4}}\beta$ is $L^2$ on $X$ and $\bar\partial\beta = \bar\partial\tilde\alpha$ on $X$.

Now $\tilde\alpha - \beta$ is a $\bar\partial$-closed $(0,1)$-form on $X$. By the vanishing of $H^1(X, \mathcal{O}_X)$ due to the simply connectedness of $X$, we can solve the equation

$$\bar\partial g = \tilde\alpha - \beta$$

to get some $L^2$ function $g$ on $X$ so that $g$ is $C^{m-1}$ on $X - M$. Take a $C^m$ family of $C^m$ real hypersurfaces $M_\sigma$ in $X$ for $-\varepsilon < \sigma < \varepsilon$ such that $M = M_0$ and there is a $C^m$ family of $C^m$ diffeomorphisms $\varphi_\sigma : M \to M_\sigma$ parametrized by $\sigma$ with $\varphi_0$ equal to the identity map of $M_0$.

We claim that $g|_{M_\sigma}$ approaches some $C^p$ function $\tilde v$ on $M$ in the $L^1$ sense as $\sigma \to 0$ (when $M_0$ is identified with $M_\sigma$ via $\varphi_\sigma$) and $\bar\partial_b \tilde v = \tilde\alpha$ on $M$. We need only verify the claim locally. We take a point $P$ in $M$ and an open neighborhood $U$ of $P$ in $X$ so that $\tilde{\alpha} - \bar{\partial}h$ vanishes to order at least $m-2$ at points of $M \cap U$ for some $C^m$ function $h$ on $U$. Consider the equation

$$\bar{\partial}(g - h) = (\tilde{\alpha} - \bar{\partial}h) - \beta$$

on $U$. The $(0,1)$-form

$$\frac{1}{\operatorname{dist}_M^{m-4}}\left((\tilde{\alpha} - \bar{\partial}h) - \beta\right)$$

is $L^2$ on $U$. By Lemma (4.1)(4), when $M_0$ is identified with $M_\sigma$ via $\varphi_\sigma$, the function $(g - h)|_{U \cap M_\sigma}$ approaches some $C^p$ function $v_U$ on $U \cap M$ in the $L^1$ sense as $\sigma \to 0$ and $\bar{\partial}_b v_U = 0$ on $U$. Since $h$ is $C^m$ on $U$ and $\tilde{\alpha} - \bar{\partial}h$ vanishes to order at least $m-2$ at points of $M \cap U$, it follows that when $M_0$ is identified with $M_\sigma$ via $\varphi_\sigma$, the function $g|_{M_\sigma}$ approaches some $C^p$ function $\tilde{v}$ on $M$ in the $L^1$ sense and $\bar{\partial}_b \tilde{v} = \tilde{\alpha}$ on $M$ as $\sigma \to 0$. Now we need only set $v = \tilde{v} - H$ on $M$. □

PROPOSITION 5.5 (analog in Levi-flat hypersurfaces of Kodaira's $\partial\bar{\partial}$ exactness). *Let $X$ be an irreducible compact Hermitian symmetric manifold of complex dimension $n$ and let $s$ be the degree of the nondegeneracy of the bisectional curvature of $X$ such that $n - s \geq 2$. Let $m, p$ be integers such that $m > \frac{3n}{2} + 4$ and $1 \leq p < m - \frac{3n}{2} - 3$. Let $M$ be a $C^{m+1}$ Levi-flat real hypersurface in $X$. Let $\omega$ be a $(1,1)$-form on the holomorphic leaves of the foliation of $M$ (i.e., $\omega$ is a section of $\left(T_M^{1,0}\right)^* \wedge \left(T_M^{0,1}\right)^*$ over $M$) such that $\omega = d\theta$ as functions on $\wedge^2(T_M^{1,0} \oplus T_M^{0,1})$ for some $C^m$ 1-form $\theta$ on the holomorphic leaves of the foliation of $M$ (i.e., $\theta$ is a function on $T_M^{1,0} \oplus T_M^{0,1}$). Then there exists a $C^p$ function $\psi$ on $M$ such that $\omega = \partial_b \bar{\partial}_b \psi$ on $M$ (i.e., $\omega = \partial\bar{\partial}\psi$ in $\left(T_M^{1,0}\right)^* \wedge \left(T_M^{0,1}\right)^*$).*

*Proof.* By taking separately the real and imaginary parts of $\omega$ and $\theta$, we can assume without loss of generality that both $\omega$ and $\theta$ are real-valued. Let $\theta = \theta^{(1,0)} + \theta^{(0,1)}$ so that $\theta^{(1,0)}$ is a function on $T_M^{1,0}$ and $\theta^{(0,1)}$ is a function on $T_M^{0,1}$. From type considerations we conclude that $\bar{\partial}_b \theta^{(0,1)} = 0$. By Proposition 5.4 we can find a $C^p$ function $g$ on $M$ such that $\bar{\partial}_b g = \theta^{(0,1)}$. Hence

$$\omega = \bar{\partial}_b \theta^{(1,0)} + \partial_b \theta^{(0,1)} = \partial_b \bar{\partial}_b (g - \bar{g}).$$  □

LEMMA 5.6. *Let $m \geq 1$. Let $Y$ be a simply connected manifold and $Z$ be a $C^m$ real hypersurface of $Y$. Then there exists a $C^m$ real-valued function $g$ on $Y$ whose zero-set is $Z$ and whose differential $dg$ is nonzero at points of $Z$.*

*Proof.* We can cover $Y$ by a simple covering of open subsets $U_j$ ($j \in J$) such that for each $j \in J$ there exists a $C^m$ real-valued function $g_j$ on $U_j$ whose zero set is $Z \cap U_j$ and whose differential $dg_j$ is nonzero at points of $Z \cap U_j$. By a simple covering we mean that $\bigcap_{\nu=1}^{k} U_{j_\nu}$ is either empty or simply connected for $j_\nu \in J$. Then $\frac{g_{j_1}}{g_{j_2}}$ is a nowhere zero real-valued function on $U_{j_1} \cap U_{j_2}$. Let

$$\varepsilon_{j_1,j_2} = \frac{\frac{g_{j_1}}{g_{j_2}}}{\left|\frac{g_{j_1}}{g_{j_2}}\right|},$$

which is either identically 1 or identically $-1$, $U_{j_1} \cap U_{j_2}$. The transition functions $\varepsilon_{j_1,j_2}$ define a principal bundle over $Y$ whose structure group is the multiplicative group consisting of two elements $\{1, -1\}$. Since $Y$ is simply connected, this principal bundle is globally trivial over $Y$. Thus there exists a function $\varepsilon_j$ on $U_j$ ($j \in J$) which is either identically 1 or identically $-1$ such that $\varepsilon_{j_1,j_2} = \varepsilon_{j_1}\varepsilon_{j_2}$ on $U_{j_1} \cap U_{j_2}$ if $U_{j_1} \cap U_{j_2}$ is nonempty. By replacing $g_j$ by $\varepsilon_j g_j$ we can assume without loss of generality that $\varepsilon_{j_1,j_2} \equiv 1$ always and $\frac{g_{j_1}}{g_{j_2}}$ is always positive on $U_{j_1} \cap U_{j_2}$ if $U_{j_1} \cap U_{j_2}$ is nonempty. Let $h_{j_1,j_2} = \log \frac{g_{j_1}}{g_{j_2}}$ which is a real-valued function on $U_{j_1} \cap U_{j_2}$ if $U_{j_1} \cap U_{j_2}$ is nonempty. Then $\{h_{j_1,j_2}\}$ defines a cocycle with coefficients in $\mathbf{R}$. By partition of unity we can find a $C^m$ real-valued function $h_j$ on $U_j$ ($j \in J$) such that $h_{j_1,j_2} = h_{j_2} - h_{j_1}$ on $U_{j_1} \cap U_{j_2}$ if $U_{j_1} \cap U_{j_2}$ is nonempty. We can now define the real-valued function $g$ on $Y$ by setting $g = g_j e^{h_j}$ on $U_j$. □

PROPOSITION 5.7 (existence of metric with zero curvature). *Let $X$ be an irreducible compact Hermitian symmetric manifold of complex dimension $n$ and let $s$ be the degree of the nondegeneracy of the bisectional curvature of $X$ such that $n - s \geq 2$. Let $m \geq \frac{3n}{2} + 7$. Let $M$ be a $C^m$ Levi-flat real hypersurface in $X$. Let $N_{M,X}^{1,0}$ be the normal bundle of the foliation of $M$ which is a $\mathbf{C}$-line bundle over $M$. Then $N_{M,X}^{1,0}$ admits a $C^2$ Hermitian metric with zero curvature along the holomorphic leaves of the foliation of $M$.*

*Proof.* By Lemma 5.6 $M$ is the zero-set of some real-valued $C^m$ function $f$ on $X$ with $df$ nonzero at points of $M$. We have the short exact sequence

$$0 \to T_M^{1,0} \to T_X^{1,0}|M \to N_{M,X}^{1,0} \to 0$$

on $M$. Consider the map

$$T_X^{1,0}|M \xrightarrow{\partial f} \mathbf{C} \times M$$

from $T_X^{1,0}|M$ to $\mathbf{C} \times M$ defined by evaluation by $\partial f$. The kernel of the map is precisely $T_M^{1,0}$. Thus the map induces a $C^{m-1}$ bundle-homomorphism

$$N_{M,X}^{1,0} \xrightarrow{\approx} \mathbf{C} \times M.$$

By identifying $N_{M,X}^{1,0}$ with the orthogonal complement of $T_M^{1,0}$ in $T_X$, we obtain a $C^{m-1}$ Hermitian metric $h$ along the fibers of $N_{M,X}^{1,0}$.

Locally on $M$ we have a $C^{m-1}$ connection for $N_{M,X}$ along the $T_M^{1,0} \oplus T_M^{0,1}$ directions of $M$ so that it agrees with $\bar{\partial}$ along the $T_M^{0,1}$ directions of $M$ and is compatible with the Hermitian metric of $N_{M,X}$. By using the partition of unity we can extend these local connections to a $C^{m-1}$ connection $D$ for the line bundle $N_{M,X}^{1,0}$ on $M$. The curvature $\Omega_D$ of the connection $D$ is a $C^{m-2}$ closed 2-form on $M$. Since the line bundle $N_{M,X}^{1,0}$ is trivial on $M$, the curvature form $\Omega_D$ is $d$-exact on $M$ and is equal to the exterior differential of a $C^{m-1}$ 1-form $\theta$ on $M$. Let $\omega$ be the $(1,1)$-component of $\Omega_D$ (i.e., $\omega$ is the restriction of $\Omega_D$ to $\left(T_M^{1,0}\right)^* \wedge \left(T_M^{0,1}\right)^*$). By Proposition (5.5) we can write it as $\omega = \partial_b \bar{\partial}_b \psi$ for some $C^2$ function $\psi$ on $M$, because $m \geq \frac{3n}{2} + 7$. So $he^{-\psi}$ is a $C^2$ metric of $N_{M,X}^{1,0}$ which has zero curvature along the holomorphic leaves of the foliation of $M$. $\square$

## 6. Curvature of quotient bundle

We now recall the domination of the curvature of a vector bundle by that of its quotient bundle, given e.g., in [G, pp. 196–199].

Let $E$ be a holomorphic vector bundle of rank $r$ with a Hermitian metric over a complex manifold $X$, $E'$ be a holomorphic subbundle of rank $s$ of $E$, and $Q$ be the quotient vector bundle $E/E'$. Let $D_E$ be a complex metric connection of $E$ (i.e., a connection which agrees with $\bar{\partial}$ in the $(0,1)$-direction and is compatible with the Hermitian metric of $E$). We choose a local orthonormal frame $e_\alpha$ ($1 \leq \alpha \leq r$) for $E$ so that $e_\alpha$ ($1 \leq \alpha \leq s$) belongs to $E'$. Write $D_E e_\alpha = \sum_{\beta=1}^r \omega_\alpha^\beta e_\beta$. Since the connection is compatible with the metric, by differentiating $\langle e_\alpha, e_\beta \rangle = 0$ or $1$ we conclude that $\omega_\alpha^\beta = -\overline{\omega_\beta^\alpha}$ for $1 \leq \alpha, \beta \leq r$.

Since $E'$ is a holomorphic subbundle of $E$, it follows that $\omega_\alpha^\beta$ is of type $(1,0)$ for $1 \leq \alpha \leq s$ and $s < \beta \leq n$, because the covariant derivatives of $e_1, \cdots, e_s$ in any $(0,1)$-direction are still spanned by $e_1, \cdots, e_s$.

We identify $Q$ with the orthogonal complement of $E'$ in $E$. Let $D_Q$ denote the complex metric connection of $Q$ defined by the metric of $Q$ induced from $E$. We have $D_Q e_\alpha = \sum_{\beta=s+1}^r \omega_\alpha^\beta e_\beta$ ($s < \alpha \leq r$). The verification is as follows. From $\omega_\alpha^\beta = -\overline{\omega_\beta^\alpha}$ for $s < \alpha, \beta \leq r$ it follows that $D_Q$ is compatible with the metric of $Q$ induced from $E$. Take a local section $\sum_{\alpha=s+1}^r f_\alpha e_\alpha$ of $E$ whose image in $Q$ is a local holomorphic section of $Q$. This means that we can find functions $f_1, \cdots, f_s$ such that $\sum_{\alpha=1}^r f_\alpha e_\alpha$ is a local holomorphic section of $E$. Let $\xi$ be a tangent vector of $X$ of type $(0,1)$. We have to verify that $\left\langle D_Q \left(\sum_{\alpha=s+1}^r f_\alpha e_\alpha \right), \xi \right\rangle$ is zero. The holomorphicity of $\sum_{\alpha=1}^r f_\alpha e_\alpha$ means that

$$\sum_{\alpha=s+1}^r (d_\xi f_\alpha) e_\alpha + \sum_{\beta=1}^r \sum_{\alpha=s+1}^r f_\alpha \omega_\alpha^\beta(\xi) e_\beta = -\sum_{\alpha=1}^s (d_\xi f_\alpha) e_\alpha - \sum_{\beta=1}^r \sum_{\alpha=1}^s f_\alpha \omega_\alpha^\beta(\xi) e_\beta.$$

Equating the part spanned by $e_{s+1}, \cdots, e_r$ and using the fact that $\omega_\alpha^\gamma$ is of type $(1,0)$ for $1 \leq \alpha \leq s$ and $s < \gamma \leq r$, we obtain

$$\left\langle D_Q\left(\sum_{\alpha=s+1}^r f_\alpha e_\alpha\right), \xi\right\rangle = \sum_{\alpha=s+1}^r (d_\xi f_\alpha) e_\alpha + \sum_{\beta=s+1}^r \sum_{\alpha=s+1}^r f_\alpha \omega_\alpha^\beta(\xi) e_\beta$$
$$= -\sum_{\beta=s+1}^r \sum_{\alpha=1}^s f_\alpha \omega_\alpha^\beta(\xi) e_\beta = 0.$$

The curvature tensors $\Theta_E$ and $\Theta_Q$ of $E$ and $Q$ are respectively given by

$$(\Theta_E)_\alpha^\beta = \frac{\sqrt{-1}}{2\pi}\left(d\omega_\alpha^\beta - \sum_{\gamma=1}^r \omega_\alpha^\gamma \wedge \omega_\gamma^\beta\right) \quad (1 \leq \alpha, \beta \leq n),$$

$$(\Theta_Q)_\alpha^\beta = \frac{\sqrt{-1}}{2\pi}\left(d\omega_\alpha^\beta - \sum_{\gamma=s+1}^r \omega_\alpha^\gamma \wedge \omega_\gamma^\beta\right) \quad (s < \alpha, \beta \leq n).$$

It follows that

$$(6.0.1) \qquad (\Theta_E)_\alpha^\beta = (\Theta_Q)_\alpha^\beta - \frac{\sqrt{-1}}{2\pi} \sum_{\gamma=1}^s \omega_\gamma^\beta \wedge \overline{\omega_\gamma^\alpha} \quad (s < \alpha, \beta \leq r).$$

PROPOSITION 6.1. $\langle \Omega_E(\xi, \bar{\xi}) e, e \rangle \leq \langle \Omega_Q(\xi, \bar{\xi}) \hat{e}, \hat{e} \rangle$ for any tangent vector $\xi$ of $X$ of type $(1,0)$ and for any $e \in E$ inducing $\hat{e} \in Q$.

Proposition 6.1 follows immediately from (6.0.1) and the fact that $\omega_\alpha^\beta$ is of type $(1,0)$ for $1 \leq \alpha \leq s$ and $s < \beta \leq n$.

## 7. Contradiction from the maximum principle

We now give the last step of the proof of Theorem 1. Assume that $M$ and $X$ are from the assumptions of Theorem 1 and we are going to derive a contradiction. Let $e^{-\varphi}$ be the Hermitian metric of $N_{M,X}^{1,0}$ induced from the standard Hermitian metric of $T_X^{1,0}$. Let $\Theta_{T_X^{1,0}}$ and $\Theta_{N_{M,X}^{1,0}}$ denote respectively the curvature forms of $T_X^{1,0}$ and $N_{M,X}^{1,0}$. From Proposition (5.7) we have a $C^2$ Hermitian metric $e^{-\psi}$ of $N_{M,X}^{1,0}$ which has zero curvature along each holomorphic leaf of the foliation of $M$. Let $\Phi$ be the function $\varphi - \psi$ on $M$. Let $P_0$ be the point of $M$ where $\Phi$ achieves its maximum value. Let $z_1, \cdots, z_{n-1}$ be the local holomorphic coordinate at $P_0$ of the holomorphic leaf through $P_0$ of the foliation of $M$ such that $\frac{\partial}{\partial z_1}, \cdots, \frac{\partial}{\partial z_{n-1}}$ are mutually orthogonal unit tangent vectors of $X$ of type $(1,0)$ at $P_0$. Let $e$ be the unit tangent vector of $X$ of type $(1,0)$ at $P_0$ which is orthogonal to $T_M^{1,0}$. Since the bisectional curvature of $X$ is $(n-2)$-nondegenerate (and in particular $(n-1)$-nondegenerate), it follows from Proposition 6.1 that

$$
\begin{aligned}
0 \geq \frac{1}{4\pi^2} \sum_{j=1}^{n-1} \frac{\partial^2}{\partial z_j \partial \overline{z_j}} \Phi &= \sum_{j=1}^{n-1} \Theta_{N^{1,0}_{M,X}} \left( \frac{\partial}{\partial z_j}, \overline{\frac{\partial}{\partial z_j}} \right) \\
&\geq \sum_{j=1}^{n-1} \left\langle \Theta_{T^{1,0}_X} \left( \frac{\partial}{\partial z_j}, \overline{\frac{\partial}{\partial z_j}} \right) e, e \right\rangle > 0
\end{aligned}
$$

at $P_0$, which is a contradiction and ends the proof of Theorem 1.

## 8. Remarks on generalizations to the higher codminsional case

For the generalization of Theorem 1 to the case where $M$ is of real codimension $q$ in $X$ and the complex dimension of $M$ is everywhere $n-q$, some modifications are needed.

The more difficult one, as mentioned in the introduction, is that Lemma 4.1 has to be suitably modified, because, when the real codimension of $M$ is greater than 1, the restriction to $M$ of a local solution of the $\bar{\partial}$-equation on $X$ with an $L^2$ right-hand side is no longer $L^1$ on $M$.

One needed modification is to get a complete Kähler metric $\omega$ on $X - M$ and a function $\psi$ on $X - M$ of the same order as $-\log \operatorname{dist}_M$ near $M$ so that the sum of an appropriate number of eigenvalues of $\psi$ with respect to $\omega$ has a positive lower bound on all of $X - M$. When the real codimension of $M$ is greater than 1, in the directions normal to $M$ the eigenvalues of $-\log \operatorname{dist}_M$ with respect to the standard Kähler metric $\omega_0$ of $X$ go to $-\infty$ of order 2 as the point approaches $M$. In the good directions tangential to $M$ the eigenvalues of $-\log \operatorname{dist}_M$ with respect to $\omega_0$ have a positive lower bound as the point approaches $M$. On the other hand, the eigenvalues of $\omega$ with respect to $\omega_0$ are expected to admit positive upper and lower bounds in directions tangential to $M$ and to go to $\infty$ of order 2 in directions normal to $M$. One has the appropriate growth orders for such a modification. (A situation similar to ours but for a complex-analytic $M$ can be found in [Sch].)

Another needed modification is as follows. For our case of an $M$ of real codimension 1 the result of Hirschowitz [H] is used in the proof of Proposition (2.4)(1) so that we can modify the function $-\log \operatorname{dist}_M$ at the conjugacy points of $M$ and still have a positive lower bound for sums of an appropriate number of eigenvalues. When the real codimension of $M$ is greater than 1, another way of modification of $-\log \operatorname{dist}_M$ at the conjugacy points of $M$ has to be used. Proposition (2.4)(1) is used also to get the Stein property assumption in Theorem (1.1). However, that use is more a matter of expediency than absolute necessity.

Harvard University, Cambridge, MA
*E-mail address*: siu@math.harvard.edu